\documentclass[12pt,reqno]{amsart}
\usepackage{amsmath} 
\usepackage{amssymb}

\def\squarebox#1{\hbox to #1{\hfill\vbox to #1{\vfill}}}

\newcommand{\M}{{\mu}}
\newcommand{\p}{{\bf M_k}}
\newcommand{\tc}{{\mathcal M}}

\newcommand{\Z}{{\mathbb Z}}

\newcommand{\RR}{{\mathbb R}}

\newcommand{\Ss}{{\bf S}}

\newcommand{\C}{{\mathbb C}}
\newcommand{\N}{{\mathbb N}}
\newcommand{\T}{{\mathbb T}}

\newcommand{\w}{{\omega}}
\newcommand{\hM}{{\widetilde{M}}}
\newcommand{\ha}{{\widetilde{a}}}
\newcommand{\oM}{{\widehat{M}}}
\newcommand{\hph}{{\widetilde{\phi}}}
\newcommand{\wi}{{{\widetilde{(\phi_k)_r}}}}
\newcommand{\wa}{{\widetilde{(a)_r}}}
\newcommand{\phik}{{\widetilde{\phi_k}}}

\newcommand{\ov}{{\widehat{T}}}
\newcommand{\wideha}{{\widetilde{T}}}

\newcommand{\E}{{\bf E}}

\theoremstyle{plain}

\newtheorem{thm}{Theorem}

\newtheorem{lem}{Lemma}

\newtheorem{deff}{Definition}

\theoremstyle{definition}

\title{Multipliers and Toeplitz operators on Banach spaces of sequences}
\vspace{6cm}

\author[ V. Petkova]{ Violeta Petkova}
\address{LABAG, Universit\'e Bordeaux I, 351, Cours de la Lib\'eration, 33405 Talence, France}\email{petkova@math.u-bordeaux1.fr}

\vspace{1cm}

\numberwithin{equation}{section}
\begin{document}
\maketitle

{\bf Abstract}: In this paper, we prove that every multiplier $M$ (i.e. every bounded operator commuting whit the shift operator S) on a large class of Banach spaces of sequences on $\Z$ is associated to a function essentially bounded by $\|M\|$ on $spec(S)$. This function is holomorphic on $\overset{\circ}{spec}(S)$, if $\overset{\circ}{spec}(S)\neq \emptyset$. Moreover, we give a simple description of $spec(S)$. We also obtain similar results for Toeplitz operators on a large class of Banach spaces of sequences on $\Z^+$. 


\vspace{0.1cm} 

\vspace{0.5cm}
\section{Introduction}
\vspace{0.5cm}
Let $E\subset \C^\Z$ be a Banach space of sequences. Denote by $S: \C^\Z \longrightarrow \C^\Z$, the shift operator defined by $Sx=(x({n-1}))_{n \in \Z}$, for $x=(x(n))_{n \in \Z} \in \C^\Z$, so that $S^{-1}x=(x({n+1}))_{n \in \Z}$. Let $F(\Z)$ be the set of sequences on $\Z$, which have a finite number of non-zero coefficients and assume that $F(\Z)$ is dense in $ E$. The elements of $F(\Z)$ will be said finite sequences. We will call multiplier on $E$ every bounded operator $M$ on $E$ such that $MSa=SMa$, for every $a \in F(\Z)$. Denote by $\M(E)$ the space of multipliers on $E$. For $z\in \T=\{z \in \C\:|\:|z|=1\}$, set $\psi_z(x)=(x(n)z^n)_{n \in \Z}$, for $x=(x(n))_{n \in \Z}$. Notice that if for $z \in \T$, $\psi_z(E)\subset E$ and if for all $n \in \Z$, the map $p_n: x\longrightarrow x(n)$ is continuous from $E$ into $\C$, then from the closed graph theorem it follows that $\psi_z$ is bounded on $E$. In this paper, we deal with Banach spaces of sequences on $\Z$ satisfying only the following very natural hypothesis:\\

(H1) The set $F(\Z)$ is dense in $E$.\\

(H2) For every $n \in \Z$, $p_n$ is continuous from $E$ into $\C$. \\

(H3) We have $\psi_z(E)\subset E, \:\forall z \in \T$ and $\sup_{z \in \T}\|\psi_z\|<+\infty$.\\

It is easy to see that if $S(E)\subset E$, then by the closed graph theorem the restriction $S\vert_E$ of $S$ to $E$ is bounded from $E$ into $E$. From now we will say that $S$ (resp. $S^{-1}$) is bounded when $S(E)\subset E$ (resp. $S^{-1}(E)\subset E$). 
If $S(E)\subset E$, we denote by $spec(S)$ the spectrum of the operator $S$ with domain $E$. If $S$ is not bounded, denote by $spec(S)$ the spectrum of $\overline{S}$. Here $\overline{S}$ is the smallest extension of $S\vert_{F(\Z)}$ as a closed operator. Recall that the domain of $\overline{S}$ is
$$D(\overline{S})=\{x\in E,\:\exists (x_n)_{n\in \Z}\subset F(\Z)\:s.t. \:x_n\longrightarrow x\: and\: Sx_n\longrightarrow y\in E\}$$
and $\overline{S}x=y$.
Our aim is to prove that every multiplier on $E$ is associated to a $L^\infty$-function on $spec(S)$, which is holomorphic on $\overset{\circ}{spec(S)}$, if $\overset{\circ}{spec(S)} \neq \emptyset$. 
The multipliers on $l^2_\w(\Z)$, when the weight $\w$ is such that $S$ and $S^{-1}$ are bounded and $\overset{\circ}{spec(S)}\neq \emptyset$ were investigated by Shields in \cite{S} (see \cite{L} for the case of the spaces $l^p(\Z)$, $p\geq 1$). Gellar considered in \cite{G} a larger class of Banach spaces with Schauder basis. On the other hand, the case of spaces $l^2_\w(\Z)$ when $\overset{\circ}{spec(S)}=\emptyset$ was treated only recently by Esterle in \cite{E}. The case of weighted spaces $L_\delta^2(\RR)$ has been examined in \cite{P1}. Nevertheless, to our best knowledge it seems that this problem for general Banach space of sequences satisfying only hypothesis (H1), (H2) and (H3) has been not yet considered in the literature.\\

Let $e_k$ be the sequence such that ${e_k}(n)=0$ if $n \neq k$ and ${e_k}(k)=1$. Define for $r>0$, $$C_r:=\{z \in \C\:|\:|z|=r\}.$$ 
From now, $E$ is a Banach space of sequences on $\Z$ satisfying (H1)-(H3).
For $M \in\M(E)$, set $\widehat{M}=M(e_0)$ and for $z \in \C$, denote by $\hM(z)$ the formal Laurent series 
$$\hM(z)=\sum_{n \in \Z}\widehat{M}(n)z^n.$$ For $M \in \M(E)$, we call $\hM$ the symbol of $M$.
Given $a \in E$, set $\ha(z)=\sum_{n \in \Z}a(n)z^n$, for $z \in \C$. It is easy to see that $Ma=\widehat{M}*a$, for $a \in F(\Z)$ and that we have on the space of formal Laurent series 
$$\widetilde{Ma}(z)=\hM(z)\ha(z),\:\:\forall z \in \C,\: \forall a \in F(\Z),$$ but it is more difficult to determine when $\hM(z)$ converges. It is natural to conjecture that $\hM(z)$ is a convergent series for $z\in spec(S)$ and that the function $\hM$ is holomorphic on $\overset{\circ}{spec(S)}$, if $\overset{\circ}{spec(S)} \neq \emptyset.$ Notice that the main difficulty is to show that $\hM$ is essentially bounded on $spec(S)$ when ${spec(S)}$ is a circle. 
For a closed operator $A$ with dense domain, denote by $\rho(A)$ the spectral radius of $A$ defined by $\rho(A)= \sup\{|\lambda|\:,\:\lambda \in spec(A)\}.$ We suppose that at least one of the operators $S$ and $S^{-1}$ is bounded.
Our main result is the following.
\begin{thm} 
$1)$ If $S$ is not bounded, but $S^{-1}$ is bounded, $\rho(S)=+\infty$ and if $S$ is bounded, but $S^{-1}$ is not bounded, 
$\rho(S^{-1})=+\infty.$\\
$2)$ We have $spec(S)=\Bigl\{ \frac{1}{\rho(S^{-1})}\leq |z|\leq \rho(S)\Bigr\}$.\\
$3)$ Let $M\in \M(E)$. For $r> 0$ such that $C_r\subset spec(S)$, we have $\hM \in L^\infty(C_r)$ and 
$$|\hM(z)|\leq \|M\|,$$ a.e. on $C_r$.\\
$4)$ If $\rho(S)>\frac{1}{\rho(S^{-1})}$, $\hM$ is holomorphic on $\overset{\circ}{spec(S)}.$
\end{thm}
If $\rho(S^{-1})=+\infty$, here $\frac{1}{\rho(S^{-1})}$ denotes $0$. 

The class of Banach spaces of sequences on $\Z$ that we consider in this paper is very general. For example (H3) is satisfied for every Banach space $E$ with a norm such that $\|(x(n))_{n \in \Z}\|=\|(|x(n)|)_{n \in \Z}\|$. We will see later that our hypothesis imply that for $x \in E$, we have 
$$\lim_{k \to +\infty}\Bigl\|x-\sum_{p=0}^k\frac{1}{k+1}\Bigl ( \sum_{n=-p}^p {x(n)}\:e_n \Bigr)\Bigr\|=0,$$ but not necessary $\lim_{p \to +\infty}\|x-\sum_{n=-p}^p {x(n)}\:e_n \|$ as it has been assumed in \cite{G}. The Example 5 below shows that this situation appears and this makes obvious the generality of our considerations. 
We will give some classical examples of Banach spaces satisfying the conditions (H1), (H2) and (H3). \\

{\bf Example 1. }\\
Let $\w$ be a weight on $\Z$, i.e. $\w$ is a positive sequence on $\Z$. Set 
$$l_\w^p(\Z)=\Bigl\{(x(n))_{n \in \Z}\in \C^\Z\:|\:\sum_{n \in \Z}|x(n)|^p\w(n)^p<+\infty\Bigr\}, \:1\leq p<+ \infty $$
and $$\|x\|_{\w,p}=\Bigl(\sum_{n \in \Z}|x(n)|^p\w(n)^p\Bigr)^{\frac{1}{p}}.$$
It is easy to see that the Banach space $l_\w^p(\Z)$ satisfies our hypothesis. \\

{\bf Example 2.} \\
For every two weights $\w_1$ and $\w_2$ and $ 1\leq p <+\infty$, $1\leq q <+\infty$, the space $l_{\w_1}^p(\Z)\cap l_{\w_2}^q(\Z)$ with the norm $\|x\|=\max\{\|x\|_{\w_1,p},\:\|x\|_{\w_2,q}\}$ satisfies also our conditions.\\

{\bf Example 3.}\\
Let ${\mathcal K}$ be a convex, non-decreasing, continuous function on $\RR^+$ such that ${\mathcal K}(0)=0$ and ${\mathcal K}(x)>0$, for $x >0.$ For example, $\mathcal {K}$ may be $x^p$, for 
$\:1\leq p< +\infty$ or $x ^{p+sin(\log(-\log( x))},\: p>1+\sqrt{2}$. Let $\w$ be a weight on $\Z$. Set
$$l_{{\mathcal K},\w}(\Z)=\Bigl\{(x(n))_{n \in \Z}\in \C^\Z \:|\:\sum_{n \in \Z}{\mathcal K}\Bigl(\frac{|x(n)|}{t}\Bigr)\w(n)<+\infty, \:\rm {for\: some\: t>0}\Bigr\}$$
and
$$\|x\|=\inf\Bigl\{ t>0 \:|\:\sum_{n \in \Z}{\mathcal K}\Bigl(\frac{|x(n)|}{t}\Bigr)\w(n)\leq 1\Bigr\}.$$
The space $l_{{\mathcal K},\w}(\Z)$, called a weighted Orlicz space (see \cite{R}, \cite{I}), is a Banach space satisfying our hypothesis. 
We can apply Theorem 1 to the multipliers on $l_{{\mathcal K},\w}(\Z)$ as well as to the spectrum of the shift on $l_{{\mathcal K},\w}(\Z)$. It seems that in the literature there are no complete results concerning the spectrum of the shift on $l_{{\mathcal K},\w}(\Z)$.
\\

{\bf Example 4.}\\
Let $(q(n))_{n \in \Z}$ be a real sequence such that $q(n) \geq 1$, for all $n \in \Z$. For $a =(a(n))_{n \in \Z}\in \C^\Z$, set 
$$\|a\|_{\{q\}}=\inf\Bigl\{ t>0\:|\:\sum_{n\in \Z}\Bigl|\frac{a(n)}{t}\Bigr|^{q(n)}\leq 1\Bigr\}.$$ 
Consider the space $l^{\{q\}}=\{a\in \C^\Z\:|\:\|a\|_{\{q\}}<+\infty\},$ which is a Banach space (see \cite{Ed}) satisfying our hypothesis. Notice that if $\lim_{n \to +\infty}|q({n+1})-q(n)|\neq 0$ and if $\sup_{n \in \Z} q(n)<+\infty$, then either $S$ or $S^{-1}$ is not bounded (see \cite{N}). \\

{\bf Example 5.}\\
Denote by $C_{[0,2\pi]}$ the space of continuous, $2\pi$-periodic, complex-valued functions on $\RR$. For $f\in C_{[0,2\pi]}$, we denote by $\hat{f}$ the sequence of Fourier coefficients of $f$. Set ${\mathcal C}=\{\hat{f}\:|\:f \in C_{[0,2\pi]}\}$ and $\|\hat{f}\|=\|f\|_\infty,$ for $f \in C_{[0,2\pi]}$. It is easy to check that the hypothesis (H1) and (H2) are satisfied by ${\mathcal C}$. For $ \alpha \in \RR$ and $f\in C_{[0,2\pi]}$, $\psi_{{e^{i\alpha}}}(\hat{f})$ is the sequence of Fourier coefficients of the function $t \longrightarrow f(t+\alpha)$. So it is clear that (H3) is satisfied by $\mathcal{C}$. Notice that in $\mathcal{C}$, $\hat{f}$ is not the limit of $\sum_{|n|\leq k}\hat{f}(n)e_n$ as $k \to +\infty$ and the space $\mathcal{C}$ is not included in the class of Banach spaces treated in \cite{G}. \\

{\bf Remark 1.}
If both $S$ and $S^{-1}$ are unbounded then Theorem 1 is not valid in general. For example, if $E=l_\w^2(\Z)$, where $\w(2n)=1$ and $\w(2n+1)=|n|+1$, for $n \in \Z$, $S$ and $S^{-1}$ are not bounded. It is easy to see that $spec(S)=\C$ and $S^2\in \M(E)$, but $\widetilde{S^2}(z)=z^2$ is obviously not bounded on $\C$. \\

In Section 3, we investigate Toeplitz operators on a general Banach space of sequences on $\Z^+ = \N$, which will be defined precisely in Definition 2 below. There are many similarities between multipliers and Toeplitz operators. We are motivated by the recent results in \cite{E} about Toeplitz operators on $l_\w^2(\Z^+)$, where $\w$ is a weight on $\Z^+$ and the results of the author (see \cite{P2}) concerning Wiener-Hopf operators on weighted spaces $L_\delta^2(\RR^+)$. Let $\E\subset \C^{\Z^+}$ be a Banach space.
Let $F(\Z^+)$ (resp. $F(\Z^-)$) be the space of the sequences on $\Z^+$ (resp. $\Z^-$) which have a finite number of non-zero coefficients. By convention, we will say that $x\in F(\Z)$ is a sequence of $F(\Z^+)$ (resp. $F(\Z^-)$) if $x(n)=0$, for $n<0$ (resp. $n>0$). 
We will assume that $\E$ is satisfying the following hypothesis:\\

(${\mathcal H}$1) The set ${F(\Z^+)}$ is dense in $ \E$.\\

(${\mathcal H}$2) For every $n \in \Z^+$, the application $p_n:x\longrightarrow x(n)$ is continuous from $\E$ into $\C$. \\

(${\mathcal H}$3) For $x=(x(n))_{n \in \Z+}\in \E$, we have $\gamma_z(x)=(z^n x(n))_{n \in \Z^+}\in \E$, for every $z \in \T$ and $\sup_{z\in\T}\|\gamma_z\|<+\infty.$\\

Notice again that if $\gamma_z(x)\in \E$, for every $x \in \E$, then $\gamma_z:\E\longrightarrow \E$ is bounded.

\begin{deff}
We define on $\C^{\Z^+}$ the operators $\Ss_1$ and $\Ss_{-1}$ as follows.
$$ For\: u\in \C^{\Z^+},\:(\Ss_1(u))(n)=0,\:if \: n = 0\:and \:(\Ss_1(u))(n)=u(n-1),\:if \:n\geq 1$$
$$(\Ss_{-1}(u))(n)=u(n+1),\: for\: n \geq 0.$$

\end{deff}
For simplicity, we note $\Ss$ instead of $\Ss_1$. Remark that we have $\Ss_{-1}\Ss=I$, however we do not have $\Ss\Ss_{-1}=I$ and this is the main technical difficulty in the analysis of the case of Toeplitz operators. It is easy to see that if $\Ss(\E)\subset \E$, then by the closed graph theorem the restriction $\Ss\vert_{\E}$ of $\Ss$ to $\E$ is bounded from $\E$ into $\E$. We will say that $\Ss$ (resp. $\Ss_{-1}$) is bounded when $\Ss(\E)\subset\E$ (resp $\Ss_{-1}(\E)\subset \E$). 
 Next, if $\Ss\vert_{\E}$ (resp. $\Ss_{-1}\vert_{\E})$ is bounded, $spec(\Ss)$ (resp. $spec(\Ss_{-1})$) denotes the spectrum of $\Ss\vert_{\E}$ (resp. $\Ss_{-1}\vert_{\E})$. 
 If $\Ss$ (resp. $\Ss_{-1}$) is not bounded, $spec(\Ss)$ (resp. $spec(\Ss_{-1})$) denotes the spectrum of the smallest closed extension of  $\Ss\vert_{F(\Z^+)}$ (resp. $\Ss_{-1}\vert_{F(\Z^+)}$).

\begin{deff}
A bounded operator on $\E$ is called a Toeplitz operator, if we have:
$$(\Ss_{-1}T\Ss)u=Tu, \: \forall u \in F(\Z^+).$$
\end{deff}
For $u \in l^2(\Z^-)\oplus \E $ introduce
$$(P^+(u))(n)=u(n),\:\forall n\geq 0\:{\rm and}\: (P^+(u))(n)=0,\:\forall n<0.$$
Given a Toeplitz operator T, set $\widehat{T}(n)=<Te_0,e_{-n}>$ and $\widehat{T}({-n})=<Te_n,e_0>$, for $n \geq 0$. Define $\widehat{T}=(\widehat{T}(n))_{n \in \Z}.$ It is easy to check that 
$$ Tu=P^+(\widehat{T}*u),\:\forall u \in F(\Z^+).$$
Set 
$$\wideha(z)=\sum_{n \in \Z}\widehat{T}(n)z^n,$$
for $z \in \C$. Notice that the series $\tilde{T}(z)$ could diverge.\\

Taking into account the similarities between multipliers and Toeplitz operators, it is natural to obtain analogous results for Toeplitz operators and to conjecture that $\wideha(z)$ converges for $z \in spec(\Ss)\cap (spec(\Ss_{-1}))^{-1}$. It is clear that if $M$ is a multiplier on ${\E}^-\oplus \E$, where $\E^-$ and $\E$ are Banach spaces of sequences respectively on $\Z^-$ and $\Z^+$, then $P^+M$ is a Toeplitz operator on $\E$. However, despite the extensive literature related to Toeplitz operators, it seems that it is not known if every Toeplitz operator is induced by a multiplier on some suitable Banach space of sequences on $\Z$. Thus we cannot use our results for the multipliers on spaces of sequences on $\Z$ to prove similar ones for Toeplitz operators. In this way, we apply the methods of Section 2 and we obtain the following theorem, when at least one of the operators $\Ss$ and $\Ss_{-1}$ is bounded. 
\begin{thm}
Let $T$ be a Toeplitz operator on $\E$.\\ 
$1)$ For $r \in \Bigl[\frac{1}{\rho({\Ss}_{-1})}, \rho(\Ss)\Bigr]$, if $\rho(\Ss)<+\infty$ or for $r \in \Bigl[\frac{1}{\rho({\Ss}_{-1})}, +\infty\Bigr[$, if $\rho(\Ss)=+\infty$ we have $\widetilde{T}\in L^\infty(C_r)$ and $|\widetilde{T}(z)|\leq \|T\|$, a.e. on $C_r$.\\
$2)$ If $\Ss$ and $\Ss_{-1}$ are bounded and if $\frac{1}{\rho({ \Ss}_{-1})}<\rho(\Ss)$, then we have $\wideha \in\mathcal{H}^\infty(\overset{\circ}{\bf {\Omega}})$, where ${\bf \Omega}:=\Bigl\{ z \in \C\:|\:\frac{1}{\rho({\Ss}_{-1})}\leq |z|\leq \rho(\Ss)\Bigr\}$.\\
$3)$ If $\Ss$ is not bounded, but ${\Ss}_{-1}$ is bounded, $\wideha \in \mathcal{H}^\infty(\overset{\circ}{U})$, where $U:=\Bigl\{ z \in \C\:|\:\frac{1}{\rho({\Ss}_{-1})}\leq |z|\Bigr\}.$
\\
$4)$ If $\Ss$ is bounded, but ${\Ss}_{-1}$ is not bounded, $\wideha \in \mathcal{H}^\infty(\overset{\circ}{V}),$ where $V:=\Bigl\{z\in \C\:|\: |z|\leq \rho(\Ss)\Bigr\}.$\\ 

\end{thm}

\vspace{0.5cm}
\section{Multipliers}
\vspace{0.5cm}
In this section, we prove Theorem 1. We denote by $E^*$ the dual space of $E$, by $\|\:.\: \| $ the norm of $E$ and by $\|\:.\:\|_* $ the norm of $E^*$. For $y \in E^*$ and $x \in E$, define $<x,y>:=y(x)$. For $k \in \Z$, setting $<x, e_k>=x(-k)$, we will consider $e_k$ as an element of $ E^*$. Notice that if we set $|||x|||=\sup_{z\in \T}\|\psi_z(x)\|$, the norm $|||.|||$ is equivalent to the norm $\|.\|$. We have $\sup\{|||\psi_z(x)|||, \:x\in E,\: |||x|||=1\}=1$, so without losing generality we can assume that $\psi_z$ is an isometry from E into E, for every $z\in \T$.
We start with the following lemma.
\begin{lem} For $x \in E$, we have 
$$\lim_{k \to +\infty}\Bigl\|\sum_{p=0}^k\frac{1}{k+1}\Bigl ( \sum_{n=-p}^p {x(n)}\:e_n \Bigr)-x\Bigr\|=0.$$
\end{lem}
{\bf Proof.}
Fix $x\in E$. The function $\Psi: z\longrightarrow \psi_z(x)$ is continuous from $\T$ into $E$. If $x \in F(\Z) $, this is obvious and if $x\in E$, the continuity follows immediately from (H1) and (H3). Consider the Fejer kernels $(g_k)_{k \in \N}\subset L^1(\T)$ defined by the formula 
$$g_k(e^{it}):=\sum_{p=0}^k\frac{1}{k+1}\sum_{|m|\leq p}e^{imt}$$
$$=\frac{1}{k+1}\Bigr(\frac{\sin(\frac{(k+1)t}{2})}{\sin\frac {t}{2}}\Bigr)^2,\:{\rm for}\: t \in \RR.$$ We have $\|g_k\|_{L^1(\T)}=1$, for $k \in \N$ and $\lim_{k \to +\infty}\int_{\delta \leq |t|\leq \pi}g_k(e^{it}) dt=0$, for $\delta >0$. Moreover, for $|n|\leq k$,
$$\hat{g_k}(n)=\frac{1}{2\pi}\int_{-\pi}^{\pi}g_k(e^{it})e^{-int}\:dt=1-\frac{|n|}{k+1}$$
and for $|n|>k,$
$$\hat{g_k}(n)=0.$$
We have $$\lim_{k \to +\infty}\|(g_k*\Psi)(1)-\Psi(1)\|=0.$$
Below we write $dz$ instead of $dm(z)$, where $m$ is the Haar measure on $\T$ such that $m(\T)=1$. For $n \in \Z$, we have
$$ \Bigl((g_k*\Psi)(1)\Bigr)(n)=\Bigl(\int_\T g_k(z)\psi_{z^{-1}}(x)dz\Bigr)(n)=\int_\T g_k(z)z^{-n}x(n) dz=\hat{g_k}(n)x(n).$$
So we obtain $$(g_k*\Psi)(1)=\sum_{n=-k}^k\Bigl( 1-\frac{|n|}{k+1}\Bigr)x(n)e_n=\sum_{p=0}^k\frac{1}{k+1}\Bigl ( \sum_{n=-p}^p {x(n)}\:e_n \Bigr)$$ and since $\Psi(1)=x$, the proof is complete.
$\Box$
\begin{lem}
For $x\in E$ and $M\in\M(E)$, the function $\tc_{x}:\T \longrightarrow E$, defined by $$\tc_{x}(z)=(\psi_{z}\circ M \circ \psi_{z^{-1}})(x)$$ is continuous.
\end{lem}
{\bf Proof.}
Fix $x$ in $F(\Z)$ and $M\in \M(E)$.
It is easy to see that
\begin{equation}\label{eq:f}
\tc_{x}(z)=(\psi_z\circ M \circ \psi_{z^{-1}})(x)=\psi_{z}(\oM)*x, \:\forall z \in \T.\end{equation} 
Indeed, for some $k\in \N$, we have
$$\Bigl((\psi_z\circ M \circ \psi_{z^{-1}})(x)\Bigr)(n)=z^n\sum_{|p|\leq k}\widehat{M}(n-p)z^{-p}x(p),\:\forall n \in \Z. $$
Thus, for every $x \in F(\Z)$, the function $z \longrightarrow (\psi_z\circ M \circ \psi_{z^{-1}})(x)$ is continuous from $\T$ into $E$. Since ${F(\Z)}$ is dense in $E$ and $\|\psi_z\circ M \circ \psi_{z^{-1}}\|\leq \|M\|$, for $z \in \T$, we deduce that $\tc_{x}$ is continuous from $\T$ into $E$, for every $x\in E$. $\Box$\\

Denote by $M_\phi$ the operator of convolution with $\phi \in F(\Z)$, when $\phi*E\subset E$. Then it is clear that $\widehat{M_\phi}=\phi$. We need the following.

\begin{lem} Let $ M\in\M(E)$, $x\in E$.\\
$1)$ We have $$\lim_{k \to +\infty}\|\p x-Mx\|=0,$$ where for $k\in \N$,
$$\p=\sum_{p=0}^k\frac{1}{k+1}\Bigl ( \sum_{n=-p}^p \widehat{M}(n)\:S^n \Bigr)=\sum_{n=-k}^{k}\Bigl(1-\frac{|n|}{k+1}\Bigr) \widehat{M}(n)S^n.$$
$2)$ $\|\p\|\leq \|M\|,$ $\forall k \in \N$.\\
$3)$ If $S^{-1}$ is not bounded, but $S$ is bounded, $\widehat{M}(n)=0,$ for $n<0$, while if $S^{-1}$ is bounded, but $S$ is not bounded, $\widehat{M}(n)=0,$ for $n>0$.

\end{lem}
{\bf Proof.} It is immediate to see that $\Bigl(1-\frac{|n|}{k+1}\Bigr) \widehat{M}(n)$ converges to $\widehat{M}(n)$, for $n \in\Z$ and assertion 1) follows from the density of $F(\Z)$ in $E$. However the control of the norm of $\p$ is less obvious. The proof follows with some modifications the arguments of \cite{S} in our more general case. Consider the Fejer kernels $(g_k)_{k \in \N}\subset L^1(\T)$ defined in the proof of Lemma 1.
Fix $ M\in\M(E)$. 
Since $\tc_{x}$ is continuous from $\T$ into $E$, for every $x \in E$ and $\tc_{x}(1)=Mx$, we have
$$\lim_{k \to +\infty}\|(g_k*\tc_{x})(1)-Mx\|=0,\:\forall x \in E.$$
Fix $x \in F(\Z).$ 
For $k \in \N$, we obtain
$$(g_k*\tc_{x})(1)=\int_{\T}g_k(z)\tc_{x}(z^{-1})dz$$
$$=\int_{\T}g_k(z)\psi_{z^{-1}}(M\psi_z (x))dz
=\int_\T g_k(z)(\psi_{z^{-1}}(\oM)*x)dz,$$
taking into account (\ref{eq:f}).
Then we have 
$$(g_k*\tc_{x})(1)=\Bigl(\int_{\T}g_k(z)\psi_{z^{-1}}(\oM)dz\Bigr)* x.$$
We observe that, for $|n|\leq k$, we have
$$\Big( \int_\T g_k(z) \psi_{z^{-1}}(\oM) dz \Big) (n)=\int_\T g_k(z)z^{-n}\widehat{M}(n)dz=
\widehat{g_k}(n){\widehat M}(n)=\Bigl(1-\frac{|n|}{k+1}\Bigr) \widehat{M}(n),$$
while for $|n|>k$, we get
$$\Big( \int_\T g_k(z) \psi_{z^{-1}}(\oM) dz \Big) (n)=0.$$ 
Since 
$$\widehat{\p}=\sum_{n=-k}^{k}\Bigl(1-\frac{|n|}{k+1}\Bigr) \widehat{M}(n)e_n,$$
it follows that
$$\widehat{\p}=\Big( \int_\T g_k(z) \psi_{z^{-1}}(\oM) dz \Big).$$
Now it is clear that

$$\|\p a\| =\|\widehat{\p}*a\|=\Bigl\|\int_\T g_k(z)(\psi_{z^{-1}}(\oM)*a)dz\Bigr\|=\Bigl\|\int_{\T}g_k(z)(\psi_{z^{-1}}\circ M\circ\psi_z)(a)dz\Bigr\|$$
$$\leq \int_\T |g_k(z)| \:\|\psi_{z^{-1}}\|\|M\|\|\psi_{z}\|\: \|a\|dz \leq \|M\|\|a\|,\:\forall a\in F(\Z) $$ 
and, since ${F(\Z)}$ is dense in $E$, we obtain $\|\p\| \leq \|M\|, \: \forall k \in \N.$\\

Suppose that $S^{-1}$ is not bounded, but $S$ is bounded. Fix $k \in \N.$ Since 
$$\p=\sum_{n=-k}^{k}\Bigl(1-\frac{|n|}{k+1}\Bigr)\widehat{M}(n)S^n$$ is bounded, the operator $S^{k-1}\p$ is bounded. We have the equality
$$ S^{k-1}\p= \Bigl(1-\frac{k}{k+1}\Bigr)\widehat{M}({-k})S^{-1}+\sum_{n=-k+1}^{k}\Bigl(1-\frac{|n|}{k+1}\Bigr)\widehat{M}(n)S^{n+k-1}$$
and taking into account that the operator $\sum_{n=-k+1}^{k}\Bigl(1-\frac{|n|}{k+1}\Bigr)\widehat{M}(n)S^{n+k-1}$ is bounded and that $S^{-1}$ is not bounded, it is clear that $\widehat{M}({-k})=0.$ In the same way, composing $\p$ and $S^p$, for $p=k-2,\:k-3,\:...., 1$, we obtain that $\widehat{M}({-n})=0$, for $n>0$. 
We can use the same argument if $S^{-1}$ is bounded but $S$ is not bounded. Thus the proof is complete.
$\Box$\\


\begin{lem} Let $\phi\in F(\Z)$ be such that $\phi*E\subset E.$\\
$1)$ If $S$ and $S^{-1}$ are bounded, then
$$|\hph(z)|\leq \|M_\phi\|,\:\forall z \in \Omega:=\Bigl\{ \frac{1}{\rho(S^{-1})}\leq |z|\leq \rho(S)\Bigr\}.$$
$2)$ If $S$ is not bounded, but $S^{-1}$ is bounded and $\phi \in F(\Z^{-})$, then
$$|\hph(z)|\leq \|M_\phi\|,\:\forall z\in \mathcal{O}:=\Bigl\{z \in \C\:|\:|z|\geq \frac{1}{\rho(S^{-1})}\Bigr\} .$$
$3)$ If $S$ is bounded, but $S^{-1}$ is not bounded and $\phi \in F(\Z^+)$, we have
$$|\hph(z)|\leq \|M_\phi\|,\:\forall z\in W:=\Bigl\{z \in \C\:|\:|z|\leq \rho(S)\Bigr\} .$$
\end{lem}
{\bf Proof.} Suppose that $S$ and $S^{-1}$ are bounded.
For $z\in spec (S)$, we have three cases:\\
Case 1.
The operator $S-zI$ is not injective. Then there exists $x \in E\backslash \{0\}$ such that $Sx=zx$.\\
Case 2.
The operator $S^*-zI$ is injective. Then the range of $S-zI$ is dense in $E$ and it is not closed. Consequently, there exists a sequence $(f_p)_{p \in \N} \subset E$ such that 
$$\lim_{p \to +\infty}\Bigl\|(S-zI) \frac {f_p}{\|f_p\|}\Bigr\|=0.$$
Case 3. The operator $S^*-zI$ is not injective. Then there exists $y \in E^*\backslash \{0\}$ such that $S^*y=zy.$ \\

Fix $z \in spec(S)$. First, assume that there exists $(h_p)_{p \in \N}\subset E$ such that 
$$\lim_{p \to +\infty}\|S h_p-z h_p\|=0\:\:{\rm and}\:\: \|h_p\|=1, \: \forall p \in \N.$$ It follows immediately that 
$$\lim_{p \to +\infty}\|S^kh_p-z^kh_p\|=0,\:\forall k \in \Z.$$
Then for $\phi\in F(\Z)$, we have for some $N>0$,
$$\|\phi*h_p-\hph(z)h_p\| \leq \sum_{k=-N}^{N}(\sup_{|k|\leq N}|\phi(k)|) \|S^k h_p-z^kh_p\|$$ and we obtain
$$\lim_{p \to +\infty}\|\phi*h_p-\hph(z)h_p\|=0.$$
Since $$|\hph(z)|=\|\hph(z)h_p\|=\|\hph(z)h_p-\phi*h_p\|+\|M_\phi h_p\|,$$ it follows that $|\hph(z)|\leq \|M_\phi\|.$\\

Now assume that there exists $y \in E^*\backslash \{0\} $ such that $S^*y=zy$. We obtain in the same way 
$$|\hph(z)|\leq \|M^*_\phi\|=\|M_\phi\|$$
and we conclude that for $\phi\in F(\Z)$, we have
$$|\hph(z)|\leq \|M_\phi\|,\:\forall z \in spec(S).$$
If $S$ is bounded and $S^{-1}$ is not bounded, the proof is similar. If $S$ is not bounded and $S^{-1}$ is bounded, we use the spectrum of $S^{-1}$ and the same arguments. Thus in the case when $\frac{1}{\rho(S^{-1})}=\rho(S)$ the proof is complete.\\

Suppose again that $S$ and $S^{-1}$ are bounded and that
$\frac{1}{\rho(S^{-1})}<\rho(S).$ Fix $\phi \in F(\Z)$.
Let $R_1>0,\:R_2>0$ be such that $R_1<R_2$ and such that the circles $C_{R_1}$ and $C_{R_2}$ with radius respectively $R_1$ and $R_2$ are included in $spec(S)$. Since $\hph$ is holomorphic on $\C\backslash\{0\}$ and $|\hph(z)|\leq \|M_{\phi}\|$, for $z \in C_{R_1}\cup C_{R_2}$, by the maximum modulus theorem we obtain 
$$|\hph(z)|\leq \|M_{\phi}\|, \:\forall z \in \Omega_{R_1,R_2}:=\Bigl\{z \in \C\:|\:R_1 \leq |z|\leq R_2\Bigr\}.$$ The inclusions $ C_{\rho(S)}\subset spec(S)$ and $C_{\frac{1}{\rho(S^{-1})}}\subset spec(S)$ imply
$$|\hph(z)| \leq \|M_\phi\|,\:{\rm for}\:z \in \Omega.$$
We complete the proof of 2) and 3) with similar arguments taking into account that if $\phi \in F(\Z^-)$, the function $z \longrightarrow \hph(z^{-1})$ is holomorphic on $\C$ and if $\phi \in F(\Z^+)$, $\hph$ is holomorphic on $\C$.
$\:\:\:\Box$\\

{\bf Proof of Theorem 1.}
Suppose that $S$ and $S^{-1}$ are bounded and let $M \in \M(E)$. Let $(\p)_{k \in \N}$ be the sequence constructed in Lemma 2 so that 
\begin{equation}\label{eq:h}\lim_{k \to +\infty}\|\p x-Mx\|=0 ,\: \forall x \in E \end{equation}
and $\|\p\|\leq \|M\|, \:\forall k \in \N.$ Set $\phi_k=\widehat{\p},$ for $k \in \N$, so that $\p=M_{\phi_k}.$ 
For $r>0$ and $a=(a(n))_{n \in \Z}\in E $, denote 
$(a)_r(n)=a(n)r^n$. Fix $r \in [\frac{1}{\rho(S^{-1})},\rho(S)]$. 
We have
$$|\wi(z)|\leq \|M_{\phi_k}\|\leq \|M\|, \:\forall z \in \T,\:\forall k \in \N.$$
We can extract from $\Bigl(\wi\Bigr)_{k \in \N}$ a subsequence which converges with respect to the weak topology $\sigma (L^{\infty}(\T), L^1(\T)) $ to a function $\nu_r \in L^\infty(\T).$ For simplicity, this subsequence will be denoted also by $\Big(\wi\Big)_{k \in \N}$. We obtain
$$\lim_{k \to +\infty}\int_\T \Big(\wi(z)g(z)-\nu_r(z)g(z)\Bigl)dz=0, \:\forall g \in L^1(\T)$$
and $\|\nu_r\|_{\infty}\leq \|M\|.$
It is clear that 
$$\lim_{k \to +\infty}\int_\T \Big(\wi(z)\wa(z)g(z)-\nu_r(z)\wa(z)g(z)\Bigl)dz=0, \:\forall g \in L^2(\T),\:\: \forall a \in F(\Z).$$
We conclude that, for $a \in F(\Z)$, $\Bigl( \wi\wa\Bigr)_{k \in \N}$ converges with respect to the weak topology of $L^2(\T)$ to $\nu_r \wa.$
Set $\widehat{\nu_r}(n)=\frac{1}{2\pi}\int_{-\pi}^{\pi}\nu_r(e^{it})e^{-itn}dt$, for $n \in \Z$ and let $\widehat{\nu_r}=(\widehat{\nu_r}(n))_{n \in \Z}$ be the sequence of the Fourier coefficients of $\nu_r$.
The Fourier transform from $l^2(\Z)$ to $L^2(\T)$ defined by
$$\mathcal{F}:l^2(\Z)\ni (a(n))_{n \in \Z} \longrightarrow \tilde{a}\vert_\T \in L^2(\T)$$
is unitary, so the sequence $\Bigl( (M_{\phi_k}a)_r \Bigr)_{k \in \N} =\Bigl((\phi_k)_r*(a)_r\Bigr)_{k \in \N}$ converges to $\widehat{\nu_r}*(a)_r$ with respect to the weak topology of $l^2(\Z)$. Taking into account (\ref{eq:h}), we obtain
$$\lim_{k \to +\infty}|<(M_{\phi_k}a)_r-(Ma)_r,\:b>|\leq\lim_{k \to +\infty}\|M_{\phi_k}a-Ma\|\:\|(b)_{r^{-1}}\|_*=0,\:\forall b \in F(\Z).$$
We deduce that 
$$(Ma)_r(n) =(\widehat{\nu_r}*(a)_r)(n),\: \forall n \in \Z,\: \forall a \in F(\Z).$$
This implies
$$(\oM)_r*(a)_r=\widehat{\nu_r}*(a)_r, \forall a \in F(\Z)$$
and 
$$(\oM)_r=\widehat{\nu_r}.$$
Thus, we have 
$$\hM(rz)=\sum_{n \in \Z}\widehat{M}(n)r^nz^n=\sum_{n \in \Z}\widehat{\nu_r}(n)z^n=\nu_r(z), \:\forall z \in \T.$$
Since $\|\nu_r\|_\infty\leq \|M\|$, it follows that the function $\hM$ is essentially bounded by $\|M\| $ on every circle included in $\Omega$. If $\rho(S)=\frac{1}{\rho(S^{-1})}$, it is clear that $spec(S)=C_{\rho(S)}=\Omega$.\\

We assume below that $\rho(S)>\frac{1}{\rho(S^{-1})}$. Since $( \phik)_{ k \in \N}$ is an uniformly bounded sequence of holomorphic functions on $\overset{\circ}{\Omega} $, we can replace $(\phik)_{ k \in \N}$ by a subsequence which converges to a function $\nu \in {\mathcal H}^{\infty}(\overset{\circ}{\Omega} )$ uniformly on every compact subset of $\overset{\circ}{\Omega}$. Thus, for $r\in ]\frac{1}{\rho(S^{-1})},\rho(S)[$,
the sequence $( \wi)_{ k \in \N}$ converges uniformly on $\T$ to the function $z \longrightarrow \nu(rz)$ and we obtain $\nu(rz)=\nu_r(z).$ We conclude that $\nu(rz)=\hM(rz)$, for $z \in \T$ and we get
$$\nu(z)=\hM(z)= \sum_{n \in \Z}\widehat{M}(n)z^n, \:{\rm for}\: z \in\overset{\circ}{\Omega}.$$
Consequently, $\hM$ is holomorphic on ${\overset{\circ}{\Omega}}$.\\

Now we will prove that $spec(S)=\Omega$.
Let $\alpha \not\in spec(S)$. Then $(S-\alpha I)^{-1}\in \M(E)$ and for $r>0$, if $C_r \subset \Omega$, there exists $\nu_r \in L^{\infty}(\T)$ such that
$$\mathcal{F}\Bigl({((S-\alpha I)^{-1}a)_r}\Bigr)(z)=\nu_r(z)\wa(z),\:\forall z \in \T,\:\forall a \in F(\Z).$$ 
Replacing $a$ by $(S-\alpha I)a$, it follows that 
$$\wa(z)=\nu_r(z)\mathcal{F}\Bigl({((S-\alpha I)a)_r}\Bigr)(z)=\nu_r(z)(rz-\alpha)\wa(z),\:\forall z \in \T,\:\forall a \in F(\Z),$$
and we get $(rz-\alpha)\nu_r(z)=1$. Suppose that $\alpha \in C_r$ i.e. $\alpha= rz_0,\:z_0\in \T$. For $\epsilon >0$, there exists $z_\epsilon \in \T$ such that $|rz_\epsilon-rz_0|\leq \epsilon$ and $|\nu_r(z_\epsilon)|\leq \|\nu_r\|_\infty$. This implies $1 \leq \epsilon \|\nu_r\|_\infty$ and we obtain a contradiction. We deduce that $C_r \subset spec(S),\:\Omega \subset spec(S)$ and $spec(S)=\Omega$.
If we suppose that $S$ or $S^{-1}$ is not bounded, we obtain the same results by the same argument replacing $\Omega$ by $\mathcal{O}$ and $W$, where $\mathcal{O}$ and $W$ are introduced in Lemma 4. Notice that when $spec(S)=\mathcal{O}$, we deduce that $\rho(S)=+\infty$ and when $spec(S)=W$, we conclude that $\rho(S^{-1})=+\infty$.
$\Box$\\


\vspace{1cm}
\section{Toeplitz operators}

\vspace{1cm}
In this section, we prove Theorem 2. In the same way, as in the proof of Lemma 1, for $x \in \E$, we obtain
$$\lim_{k \to +\infty}\Bigl\|\sum_{n=0}^k\frac{1}{k+1}\sum_{p=0}^nx(p)e_p-x\Bigr\|=0.$$

If $\phi\in F(\Z)$ is such that $P^+(\phi*\E)\subset \E$, we denote by $T_\phi$ the operator on $\E$ defined by $T_\phi x=P^+(\phi*x),$ for $x\in \E.$
By the same method, as in Section 2, we obtain the following lemma.

\begin{lem}
$1)$ Given a Toeplitz operator $T$ on $\E$, the sequence $(\phi_n)_{n \in \N}$, where $$\phi_n =\sum_{p=0}^n\frac{1}{n+1}\Bigl ( \sum_{k=-p}^p \widehat{T}(k)e_k \Bigr)$$ 
has the properties 
$$\lim_{n \to +\infty}\|T_{\phi_n}x-Tx\|, \:\forall x \in \E,\:and \: \|T_{\phi_n}\|\leq \|T\|,\:\forall n \in \N.$$\\
$2)$ If $\Ss$ is bounded, but ${\Ss}_{-1}$ is not bounded, $\widehat{T}(k)=0$, for $k<0$.\\
$3)$ If $\Ss$ is not bounded, but ${\Ss}_{-1}$ is bounded, $\widehat{T}(k)=0$, for $k>0$.\\

\end{lem}

\begin{lem}
$1)$ If $\Ss$ and ${\Ss}_{-1}$ are bounded, for $\phi\in F(\Z)$, we have
$$|\widetilde{\phi}(z)|\leq \|T_\phi\|, \:\forall z \in {\bf \Omega}:=\Bigl\{ z \in \C\:|\:\frac{1}{\rho({\Ss}_{-1})}\leq |z|\leq \rho(\Ss)\Bigr\}.$$
$2)$ If $\Ss$ is not bounded, but ${\Ss}_{-1}$ is bounded, for $\phi\in F(\Z^-)$, we have
$$|\widetilde{\phi}(z)|\leq \|T_\phi\|, \:\forall z \in V:=\Bigl\{ z \in \C\:|\:\frac{1}{\rho({\Ss}_{-1})}\leq |z|\Bigr\}.$$ 
$3)$ If $\Ss$ is bounded, but ${\Ss}_{-1}$ is not bounded, for $\phi\in F(\Z^+)$, we have
$$|\widetilde{\phi}(z)|\leq \|T_\phi\|, \:\forall z \in U:=\Bigl\{z\in \C\:|\: |z|\leq \rho(\Ss)\Bigr\}.$$ 

\end{lem}

{\bf Proof of Lemma 6.} We will present only the proof of 1). The proofs of 2) and 3) are very similar. Suppose that $\Ss$ and ${\Ss}_{-1}$ are bounded. Let $\lambda \in spec(\Ss)\cap(spec({\Ss}_{-1}))^{-1}.$ 
Since $\lambda \in spec(\Ss)$,
there exists a sequence $(f_n)_{n \in \N}$, $f_n\in \E$ such that 
\begin{equation}\label{eq:1} \lim_{n \to +\infty}\|\Ss f_n-\lambda f_n\|=0 \: {\rm and}\: \|f_n\|=1,\: \forall n \in \N
\end{equation}
or there exists \begin{equation}\label{eq:2}a \in \E^*\backslash\{0\},\:\Ss^*a=\lambda a.\end{equation}
If (\ref{eq:1}) holds, we obtain
\begin{equation}\label{eq:3}
\lim_{n \to +\infty}\|\Ss^kf_n-\lambda^k f_n\|=0 \:{\rm and}\:\lim_{n \to +\infty}\|\Ss_{-1}^kf_n-\lambda^{-k} f_n\|=0,\:\forall k \in \N.
\end{equation}
Since $\lambda^{-1} \in spec({\Ss}^*_{-1})$, there exists a sequence $(g_n)_{n \in \N}$, $g_n\in \E^*$ such that

\begin{equation}\label{eq:4} \lim_{n \to +\infty}\|{\Ss}^*_{-1}g_n-\lambda^{-1} g_n\|_*=0 \: {\rm and}\: \|g_n\|_*=1,\: \forall n \in \N
\end{equation}
or there exists \begin{equation}\label{eq:5}b \in \E\backslash\{0\},\:({\Ss}^*_{-1})^*b=\Ss_{-1}b=\lambda^{-1} b.\end{equation}
Next if (\ref{eq:4}) holds, we get
\begin{equation}\label{eq:6}
\lim_{n \to +\infty}\|({\Ss}^*)^{k}g_n-\lambda^k g_n\|_*=0 \:{\rm and}\:\lim_{n \to +\infty}\|({\Ss}_{-1}^*)^{k}g_n-\lambda^{-k} g_n\|_*=0 , \:\forall k \in \N.
\end{equation}
Suppose that we have (\ref{eq:2}) and (\ref{eq:5}). Let $a\in \E^*\backslash\{0\}$ be such that ${\Ss}^*(a)=\lambda a$. Set
$$x(-n)=<e_n,x>=<\Ss^ne_0,x>=<e_0,\Ss^{*n}x>,\:\forall x\in E^*.$$
Since $F(\Z^+)$ is dense in $\E$, the map
$$E^*\ni x\longrightarrow (x(-n))_{n\geq 0}$$
is injective. We have $$a(-n)=\lambda^n a(0),\:\:n\geq 0.$$ 
Let $b\in \E\backslash\{0\}$ be such that $\Ss_{-1}b=\lambda^{-1} b$. We obtain 
$$b(n)\lambda^n=b(0), \:n\geq 0.$$
Since $a\neq 0$ and $b \neq 0$, we have $a(0)\neq 0$, $b(0)\neq 0$.
For $k \in \N$, define $u_k \in F(\Z^+)$ by
$$u_k=\sum_{n=0}^k\frac{1}{k+1}\sum_{p=0}^nb(p)e_p=\sum_{n=0}^k\Bigl(1-\frac {n}{k+1}\Bigr)b(n)e_n.$$
We have $\lim_{k \to +\infty}\|u_k-b\|=0$ and so $\lim_{k \to +\infty}<u_k,a>=<b,a>$. On the other hand,
$$\lim_{k \to +\infty}<u_k,a>=\lim_{k \to +\infty}\sum_{n=0}^k\Bigl(1-\frac {n}{k+1}\Bigr)\lambda^{-n}b(0)\lambda^na(0)=\lim_{k \to +\infty}\Bigl(\frac{k}{2}+1\Bigr)a(0)b(0)=+\infty.$$ We obtain an obvious contradiction and
we conclude that we cannot have in the same time (\ref{eq:2}) and (\ref{eq:5}), hence we have (\ref{eq:3}) or (\ref{eq:6}). Using the same arguments as in the proof of Lemma 4 and (\ref{eq:3}) or (\ref{eq:6}), we deduce 
$$|\widetilde{\phi}(\lambda)|\leq \|T_\phi\|,\:\forall \phi \in F(\Z), \:\forall \lambda \in spec(\Ss)\cap(spec({\Ss}_{-1}))^{-1}.$$
By the maximum modulus theorem we obtain 

\begin{equation}\label{eq:7}
|\widetilde{\phi}(\lambda)|\leq \|T_\phi\|,\:\forall \phi \in F(\Z), \:\forall \lambda \in {\bf \Omega}.
\end{equation}
If $\Ss$ is bounded and ${\Ss}_{-1}$ is not bounded, then for $\lambda \in spec(\Ss)$ there exists a sequence $(h_n)_{n \in \N}$, $h_n\in \E$ such that $\lim_{n \to +\infty}\|\Ss h_n-\lambda h_n\|=0 $ and $\|h_n\|=1$ or there exists $c \in \E^*\backslash \{0\}$ such that $\Ss^*c=\lambda c$. Using the same arguments as in the proof of Lemma 4, we obtain 
$$|\widetilde{\phi}(\lambda)|\leq \|T_\phi\|,\:\forall \phi \in F(\Z^+), \:\forall \lambda \in spec(\Ss).$$
If ${\Ss}_{-1}$ is bounded, we use the spectrum of ${\Ss}_{-1}$. In the both situations, we obtain the result by using the maximum modulus theorem. $\Box$\\

Now we will prove the main result in this section.\\

{\bf Proof of Theorem 2.} The proof of Theorem 2 goes by using the same arguments as the proof of Theorem 1 with minor modifications. For the convenience of the reader we will give the main steps.
First, assume that $\Ss$ and ${\Ss}_{-1}$ are bounded.
Let $T$ be a Toeplitz operator on $\E$ and let $(\phi_k)_{k \in \N}\subset F(\Z)$ be such that
$$\lim_{k \to +\infty} \|T_{\phi_k}a-Ta\|=0, \:\forall a \in \E$$
and
$$\|T_{\phi_k}\|\leq \|T\|, \:\forall k \in \N.$$ 
For $r>0$ and $a \in \E$, denote 
$(a)_r(n)=a(n)r^n.$ Fix $r \in [\frac{1}{\rho({\Ss}_{-1})},\rho(\Ss)]$. 
We have
$$|\widetilde{(\phi_k)_r}(z)|\leq \|T_{\phi_k}\|\leq \|T\|, \:\forall z \in \T,\:\forall k \in \N.$$
We can extract from $\Bigl(\widetilde{(\phi_k)_r}\Bigr)_{k \in \N}$ a subsequence which converges with respect to the weak topology $\sigma(L^{\infty}(\T), L^1(\T)) $ to a function $\nu_r \in L^\infty(\T).$ For simplicity, this subsequence will be denoted also by $\Big(\widetilde{(\phi_k)_r}\Big)_{k \in \N}$. 

We conclude that, for $a \in F(\Z)$, $\Bigl( \widetilde{(\phi_k)_r}\widetilde{(a)_r}\Bigr)_{k \in \N}$ converges with respect to the weak topology of $L^2(\T)$ to $\nu_r \widetilde{(a)_r}.$
Denote by $\widehat{\nu_r}=(\widehat{\nu_r}(n))_{n \in \Z}$ the sequence of the Fourier coefficients of $\nu_r$.
Since the Fourier transform from $l^2(\Z)$ to $L^2(\T)$ is an isometry, the sequence $(\phi_k)_r*(a)_r$ converges to $\widehat{\nu_r}*(a)_r$ with respect to the weak topology of $l^2(\Z)$. On the other hand, $\Big( T_{\phi_k}a \Big)_{k \in \N} $ converges to $Ta$ with respect to the topology of $\E$. Consequently, we have 
$$\lim_{k \to +\infty}|<(T_{\phi_k}a)_r-(Ta)_r,\:e_{-n}>|$$
$$\leq\lim_{k \to +\infty}\|T_{\phi_k}a-Ta\|\:\|(e_{-n})_{r^{-1}}\|_*=0,\:\forall n \in\N,\:\forall a \in F(\Z^+).$$
We conclude that 
$$(Ta)_r =P^+(\widehat{\nu_r}*(a)_r),\: \forall a \in F(\Z^+).$$
Since
$$(Ta)_r =P^+((\ov*a)_r),\: \forall a \in F(\Z^+),$$ it follows that
$\widehat{T}(n)r^n=\widehat{\nu_r}(n),\:\forall n \in \Z.$
From the estimation $\|\nu_r\|_\infty\leq \|T\|$, we deduce that the function $\wideha$ is essentially bounded by $\|T\| $ on every circle included in ${\bf \Omega}$. \\

If we assume that $\rho(\Ss)>\frac{1}{\rho({\Ss}_{-1})}$, as in the proof of Theorem 1, we conclude that $\wideha$ is holomorphic on ${\overset{\circ}{{\bf \Omega}}}$.\\

Replacing ${\bf \Omega}$ by $U$ and $V$ and using the same arguments, we obtain the results when one of the operators $\Ss$ and ${\Ss}_{-1}$ is not bounded.
$\Box$\\
\vspace{0.3cm}

{\bf Acknowledgments.} The author thanks Jean Esterle for his useful advices and encouragements.\\
\end{document}